\documentclass[12pt]{amsart} 
\usepackage[mathscr]{eucal}
\usepackage{amsmath,amsfonts}

\parskip=\smallskipamount

\newtheorem{theorem}{Theorem}[section]

\newtheorem{corollary}[theorem]{Corollary}
\newtheorem{lemma}[theorem]{Lemma}

\makeatletter
\@addtoreset{equation}{section}
\makeatother

\newcommand{\CC}{{\mathbb C}}
\newcommand{\NN}{{\mathbb N}}

\newcommand{\FF}{{\mathbb F}}

\newcommand{\cB}{{\mathcal B}}

\newcommand{\cD}{{\mathcal D}}
\newcommand{\cE}{{\mathcal E}}
\newcommand{\cF}{{\mathcal F}}
\newcommand{\cG}{{\mathcal G}}
\newcommand{\cH}{{\mathcal H}}
\newcommand{\cK}{{\mathcal K}}
\newcommand{\cL}{{\mathcal L}}
\newcommand{\cM}{{\mathcal M}}
\newcommand{\cN}{{\mathcal N}}
\newcommand{\cP}{{\mathcal P}}

\newcommand{\cS}{{\mathcal S}}
\newcommand{\cT}{{\mathcal T}}

\newcommand{\cJ}{{\mathcal J}}

\newdimen\expt
\def\boxit#1{\setbox0\hbox{$\displaystyle{#1}$}
      \hbox{\lower.4\expt
 \hbox{\lower3\expt\hbox{\lower\dp0
      \hbox{\vbox{\hrule height.4\expt
 \hbox{\vrule width.4\expt\hskip3\expt
      \vbox{\vskip3\expt\box0\vskip2\expt}%
 \hskip3\expt\vrule width.4\expt}\hrule height.4\expt}}}}}}
\begin{document}
\pagestyle{plain}

\bigskip

\title 
{Tensor algebras and displacement structure. \\
II. Noncommutative Szeg\"o polynomials} 
\author{T. Constantinescu} \author{J. L. Johnson} 

\address{Department of Mathematics \\
  University of Texas at Dallas \\
  Box 830688, Richardson, TX 75083-0688, U. S. A.}
\email{\tt tiberiu@utdallas.edu} 
\address{Department of Mathematics \\
  University of Texas at Dallas \\
  Box 830688, Richardson, TX 75083-0688, U. S. A   } 
\email{\tt jlj@utdallas.edu}

\maketitle

\noindent
{\small {\bf Abstract.}
In this paper we continue to explore the connection 
between tensor algebras and 
displacement structure. We focus on recursive orthonormalization
and we develop an analogue of the Szeg\"o type theory of orthogonal
polynomials in the unit circle for several noncommuting variables. 
Thus, we obtain the recurrence equations and Christoffel-Darboux
formulas for Szeg\"o polynomials in several noncommuting 
variables, as well as a Favard type result. 
Also we continue to study a Szeg\"o type kernel
for the $N$-dimensional unit ball of an infinite dimensional
Hilbert space.


}

\section{Introduction}

In the first part of this paper, \cite{CJ2}, 
we explored the connection between tensor algebras and displacement 
structure. 
The displacement structure theory was initiated in \cite{KKM}
as a recursive factorization theory for matrices
whose implicit structure is encoded by a so-called 
displacement equation. This has been useful in several
directions including constrained and unconstrained 
rational interpolation, maximum entropy, inverse scattering,
$H^{\infty }$-control, signal detection, digital filter
design, nonlinear Riccati equations, certain Fredholm and 
Wiener-Hopf equations, etc., 
see \cite {KS}. Aspects of the Szeg\"o theory
can be also revealed within the displacement structure
theory.  Our main goal is to develop an analogue for 
polynomials in several noncommuting variables of the Szeg\"o theory
of orthogonal polynomials on the unit circle.
An analogue of the Szeg\"o theory of orthogonal polynomials
on the real line is being developed in the companion paper
\cite{CJ3}.

The paper is organized as follows.
In Section~2 we review notation and several results from  
\cite{CJ2}. In this way, this paper can be 
read independently of \cite{CJ2}. In Section~3 we introduce
orthogonal polynomials in several noncommuting variables associated to 
certain representations of the free semigroup and discuss
their algebraic properties, mostly related to the
recursions that they satisfy. 
In Section~4 we consider several positive definite
kernels on the $N$-dimensional unit ball of an 
infinite dimensional Hilbert space. In particular, we 
prove a basic property of the Szeg\"o type 
kernel studied in \cite{CJ1}  
by characterizing its Kolmogorov decomposition.
In Section~5 we discuss the problem of recovering the 
representation from orthogonal polynomials and we prove a
Favard type result. We plan a more detailed study of 
applications to multiscale systems in a sequel of this paper.

\section{Preliminaries}
We briefly review several constructions of the tensor algebra and 
introduce the necessary notation. We also review the connection 
with displacement structure theory as established in \cite{CSK}.

\bigskip
\noindent
{\it 2.1. Tensor algebras}. The tensor algebra over $\CC ^N$ 
is defined by the algebraic direct sum 
$$\cT _N=\oplus _{k\geq 0}(\CC ^N)^{\otimes k},$$
where $(\CC ^N)^{\otimes k}$ denotes the $k$-fold tensor
product of $\CC ^N$ with itself. The addition is the 
componentwise addition and the  multiplication is defined
by juxtaposition:
$$(x\otimes y)_n=\sum _{k+l=n}x_k\otimes y_l.$$
If $\{e_1,\ldots ,e_N\}$ is the standard basis of
$\CC ^N$, then $\{e_{i_1}\otimes \ldots \otimes e_{i_k}
\mid 1\leq i_1, \ldots ,i_k\leq N\}$ is an orthonormal 
basis of $\cT _N$. 
Let $\FF _N^+$ be the unital free semigroup 
on $N$ generators $1,\ldots ,N$ with lexicograhpic
order $\prec $. The empty word is the identity element
and the length of the word $\sigma $ is 
denoted by $|\sigma |$. The length of the empty word is $0$.
If $\sigma =i_1\ldots i_k$ then we write
$e_{\sigma }$ instead of $e_{i_1}\otimes \ldots \otimes e_{i_k}$, 
so that any element of $\cT _N$ can be uniquely written in the 
form
$x=\sum _{\sigma \in \FF _N^+}c_{\sigma }e_{\sigma }$,
where only finitely many of the complex numbers $c_{\sigma }$ 
are different from $0$.

Another construction of $\cT _N$ can be obtained as follows.
Let $S$ be a unital semigroup and denote by $F_0(S)$
the set of functions $\phi :S\rightarrow \CC $ with the property
that $\phi (s)\ne 0$ for only finitely many values of $s$. 
This set has a natural vector space structure and 
$B_S=\{\delta _s \mid s\in S\}$ is a vector basis for 
$F_0(S)$, where $\delta _s$ is the Kronecker symbol associated
to $s\in S$. Also $F_0(S)$ is a unital associative algebra
with respect to the product
$$\begin{array}{rcl}
\phi *\psi &=&(\sum _{s\in S}\phi (s)\delta _s)*
(\sum _{t\in S}\psi (t)\delta _t) \\
 & & \\
 &=&\sum _{s,t\in S}\phi (s)\psi (t)\delta _{st}.
\end{array}
$$
It is readily seen that $F_0(\FF _N^+)$ is isomorphic to 
$\cT _N$. Since each element in $F_0(\FF _N^+)$
can be uniquely written as a (finite) sum,
$\phi =\sum _{\sigma \in \FF _N^+}c_{\sigma }\delta_{\sigma },$
the isomorphism is the linear extension $\Phi _1$ of the mapping 
$\delta _{\sigma }\rightarrow e_{\sigma }$, $\sigma \in \FF _N^+$.

Another copy of the tensor algebra is given by the 
algebra $\cP _N^0$ of polynomials in $N$ noncommuting
indeterminates $X_1,\ldots ,X_N$ with complex coefficients.
Each element $P\in \cP _N^0$ can be uniquely written 
in the form $P=\sum _{\sigma \in \FF _N^+}c_{\sigma }X_{\sigma }$,
with $c_{\sigma }\ne 0$ for finitely many $\sigma $'s and 
$X_{\sigma }=X_{i_1}\ldots X_{i_k}$ where $\sigma =i_1\ldots i_k
\in \FF _N^+$. The linear extension $\Phi _2$ of the mapping 
$\delta _{\sigma }\rightarrow X_{\sigma }$, $\sigma \in \FF _N^+$, 
gives an isomorphism of $\cT _N$ with $\cP _N^0$.

Yet another copy of $\cT _N$ inside the algebra of lower triangular
operators allowed for the connection with displacement structure
established in \cite{CSK}. Thus, let $\cE $ be a Hilbert space
and define: $\cE _0=\cE $ and for $k\geq 1$, 
\begin{equation}\label{prebaza}
\cE _k=\underbrace{\cE _{k-1}\oplus \ldots \oplus \cE _{k-1}}_
{N \,terms}=\cE ^{\oplus N}_{k-1}.
\end{equation}
For $\cE=\CC $ we have that $\CC _k$
can be identified with $(\CC ^N)^{\otimes k}$ and $\cT _N$ is
isomorphic to the algebra $\cL _N^0$ of lower triangular 
operators $T=[T_{ij}]\in \cL (\oplus _{k\geq 0}\CC _k)$
with the property 
\begin{equation}\label{baza}
T_{ij}=\underbrace{T_{i-1,j-1}\oplus \ldots \oplus T_{i-1,j-1}}
_{N \, terms}=T^{\oplus N}_{i-1,j-1},
\end{equation}
for $i\leq j$, $i,j\geq 1$, and $T_{j0}=0$ for all 
sufficiently large $j's$. The isomorphism is given by the 
map $\Phi _3$  defined as follows:
let $x=(x_0,x_1,\ldots )\in \cT _N$ ($x_p\in (\CC ^N)^{\otimes p}$
is the $p$th homogeneous component of $x$); then 
$x_p=\sum _{|\sigma |=p}c_{\sigma }e_{\sigma }$
and for $j\geq 0$, $T_{j0}$ denotes the column matrix
$[c_{\sigma }]^T_{|\sigma |=j}$, where $"T"$ denotes the 
matrix transpose. Then $T_{j0}=0$ for 
all sufficiently large $j$'s and we can define 
$T\in \cL (\oplus _{k\geq 0}\CC _k)$
by using \eqref{baza}. Finally, set $\Phi _3(x)=T$.

\bigskip
\noindent
{\it 2.2. Displacement structure}.
We can now describe the displacement structure of the 
tensor algebra. We write this connection for $\cL _N^0$.
Then it can be easily translated into any other
realization of the tensor algebra. Let 
$F_k=[T^k_{ij}]\in \cL (\oplus _{k\geq 0}\CC _k)$, $k=1,\ldots ,N$,
be isometries defined by the formulae: $T_{ij}=0$
for $i\ne j+1$ and $T_{i+1,i}$ is a block-column matrix
consisting of $N$ blocks of dimension dim $\CC _i$, all of them zero except
for the $k$th block which is the identity on $\CC _k$.
We have the following result noticed in \cite{CSK}.

\begin{theorem}\label{displacement}
Let $T\in \cL _N^0$ and define $A=I-TT^*$. Then
\begin{equation}\label{dispeq}
A-\sum _{k=1}^NF_kAF_k^*=GJ_{11}G^*,
\end{equation}
where 
$$G=\left[\begin{array}{cc}
1 & T_{00} \\
0 & T_{01} \\
\vdots & \vdots 
\end{array}
\right]\quad \mbox{\it and} \quad J_{11}=
\left[\begin{array}{cc}
1 & 0 \\
0 & -1
\end{array}
\right].$$
\end{theorem}

The model $\cL _N^0$ of the tensor algebra is also useful
in order to extend this algebra to some topological 
tensor algebras - see \cite{Ho1}. Here we consider 
only the norm topology and denote by $\cL _N$ the algebra
of all lower triangular operators 
$T=[T_{ij}]\in \cL (\oplus _{k\geq 0}\CC _k)$
satisfying \eqref{baza}. 

\bigskip
\noindent
{\it 2.3. Multiscale processes}.
Multiscale processes are stochastic processes indexed 
by nodes on a tree. They became quite popular lately, 
see \cite{BBCGNW}, \cite{CWB}, and have potential to model 
the self-similarity of fractional Brownian motion leading 
to iterative algorithms in computer vision, remote sensing, etc. 

Here we restrict our attention 
to the case of the Cayley tree, in which each node has $N$
branches.
The vertices of the Cayley tree are indexed by $\FF _N^+$.

Let $(X,\cF ,P)$ be a probability space and let
$\{v_{\sigma }\}_{\sigma \in \FF _N^+}\subset L^2(P)$
be a family of random variables. Its covariance kernel is
$$K(\sigma , \tau )=\int _X\overline{v}_{\sigma }v_{\tau }dP,$$
and assume that the process is stationary in the sense 
(considered earlier, e.g. \cite{Fr}) that:
\begin{equation}\label{sta1}
K(\tau \sigma ,\tau \sigma ')=K(\sigma ,\sigma '), \quad 
\tau ,\sigma ,\sigma '\in \FF _N^+,
\end{equation}
\begin{equation}\label{sta2}
K(\sigma ,\tau )=0 \quad \mbox{ if there is no $\alpha \in \FF _N^+$
such that $\sigma =\alpha \tau $ or $\tau =\alpha \sigma $}.
\end{equation}
Conversely, by the invariant Kolmogorov decomposition theorem, 
see e.g., \cite{Pa}, Ch.~II, there exists an isometric 
representation $u$ of $\FF _N^+$ on a Hilbert space $\cK $
and a mapping $v:\FF _N^+\rightarrow \cK$ such that
$$K(\sigma ,\tau )=\langle v(\tau ),v(\sigma )\rangle _{\cK },$$ 
$$u(\tau )v(\sigma )=v(\tau\sigma ),$$
for all $\sigma ,\tau \in \FF _N^+$, the set
$\{v(\sigma )\mid \sigma \in \FF _N^+\}$ is total in $\cK $, and 
$u(1)$, $\dots $, $u(N)$ are isometries with orthogonal ranges.

This class of multiscale processes would be suitable to model
branching processes without "past". If a "past"
should be attached to a process as above, we could try to
consider processes indexed by the nodes of the tree
associated to the free group on $N$ generators $1$, $\ldots $, $N$.
As mentioned in the introduction, we plan to look at this matter
in a sequel of this paper. Here we focus on processes with 
covariance kernel satisfying \eqref{sta1} and \eqref{sta2}.
It was shown in \cite{CJ2} that such a kernel has displacement
structure. Also, it is clear that for all  $j,k\geq 1$,
\begin{equation}\label{explozia}
\left[K(\sigma ,\tau )
\right]_{|\sigma |=j,|\tau |=k}=
(\left[K(\sigma ',\tau ')
\right]_{|\sigma '|=j-1,|\tau '|=k-1})^{\oplus N},
\end{equation}
so that the kernel is determined by 
the elements $s_{\sigma }=K(\emptyset ,\sigma )$, 
$\sigma \in \FF _N^+$.

By Theorem~1.5.3 in \cite{Co}, each 
positive definite kernel $K$ on $\FF _N^+$ 
is uniquely determined by a family of contractions
$\{\gamma _{\sigma , \tau }\mid \sigma ,\tau \in \FF _N^+, 
\sigma \preceq \tau \}$
such that $\gamma _{\sigma ,\sigma  }=0$, $\sigma \in \FF _N^+$,
and otherwise 
$\gamma _{\sigma, \tau }\in \cL (\cD _{\gamma _{\sigma +1,\tau}},
\cD _{\gamma ^*_{\sigma ,\tau -1}})$ (for a contraction 
$T$ between two Hilbert spaces $D_T=(I-T^*T)^{1/2}$ denotes 
the defect operator of $T$ and $\cD _T$ is the defect
space of $T$ defined as the closure of the range of $D_T$ --
note that in our case $\gamma _{\sigma , \tau }$ are just complex
numbers and the condition $\gamma _{\sigma, \tau }\in 
\cL (\cD _{\gamma _{\sigma +1,\tau}},
\cD _{\gamma ^*_{\sigma ,\tau -1}})$ for $\sigma \prec \tau $ encodes the fact
that $|\gamma _{\sigma +1,\tau}|=1$ or
$|\gamma _{\sigma ,\tau -1}|=1$ implies that $\gamma _{\sigma , \tau }=0$;
also, $\tau -1$ denotes the predecessor of $\tau $ with respect to the 
lexicographic order $\prec $ on $\FF _N^+$, while $\sigma +1$
denotes the successor of $\sigma $).
In addition, the 
positive definite kernel $K$ satisfies \eqref{sta1} and \eqref{sta2}
if and only if
\begin{equation}\label{sta3}
\gamma _{\tau \sigma ,\tau \sigma '}=\gamma _{\sigma ,\sigma '}, \quad 
\tau ,\sigma ,\sigma '\in \FF _N^+,
\end{equation}
\begin{equation}\label{sta4}
\gamma _{\sigma ,\tau }=0 \quad \mbox{ if there is no $\alpha \in \FF _N^+$
such that $\sigma =\alpha \tau $ or $\tau =\alpha \sigma $}.
\end{equation}
We define $\gamma _{\sigma }=
\gamma _{\emptyset ,\sigma }$, $\sigma \in \FF _N^+$
and we notice that 
$\{\gamma _{\sigma , \tau }\mid \sigma ,\tau \in \FF _N^+, \sigma \preceq
\tau \}$
is uniquely determined by $\{\gamma _{\sigma }\}_{\sigma \in \FF _N^+}$
by the formula
\begin{equation}\label{explo}
\left[\gamma _{\sigma ,\tau }
\right]_{|\sigma |=j,|\tau |=k}=
(\left[\gamma _{\sigma ',\tau '}
\right]_{|\sigma '|=j-1,|\tau '|=k-1})^{\oplus N}, \quad j,k\geq 1.
\end{equation}

\section{Szeg\"o polynomials}
We introduce polynomials in several noncommuting variables
orthogonal with respect to a positive definite
kernel $K$ satisfying \eqref{sta1} and \eqref{sta2}. We extend 
some elements of the Szeg\"o theory to this setting.

The kernel $K$ being given, we can introduce an inner product
on $F_0(\FF _N^+)$ in the usual manner:
\begin{equation}\label{sca}
\langle \phi ,\psi \rangle _K=\sum _{\sigma ,\tau \in \FF_N^+}
K(\sigma ,\tau )\psi (\tau )\overline{\phi (\sigma )}.
\end{equation}
By factoring out the subspace $\cN _K=\{\phi \in F_0(\FF _N^+)
\mid \langle \phi ,\phi \rangle _K=0\}$ and completing with respect
to the norm induced by \eqref{sca} we
obtain a Hilbert space denoted $\cH _K$.
A similar structure can be introduced on $\cP ^0_N$.
Let $P=\sum _{\sigma \in \FF _N^+}c_{\sigma }X_{\sigma }$, 
$Q=\sum _{\sigma \in \FF _N^+}d_{\sigma }X_{\sigma }$
be elements in $\cP ^0_N$, then define:
\begin{equation}\label{sca2}
\langle P,Q\rangle _K=\sum _{\sigma ,\tau \in \FF_N^+}
K(\sigma ,\tau )d_{\tau }\overline{c}_{\sigma }.
\end{equation}
By factoring out the subspace $\cM _K=\{P\in \cP ^0_N
\mid \langle P,P \rangle _K=0\}$ and completing with respect
to the norm induced by \eqref{sca2} we
obtain a Hilbert space denoted $L^2(K)$.
One can check that the map $\Phi _2$ defined by $\delta _{\sigma }
\rightarrow X_{\sigma }$, $\sigma \in \FF _N^+$, extends to
a unitary operator from $\cH _K$ onto $L^2(K)$.

From now on we assume that for any $\alpha \in \FF _N^+$
the matrix $\left[K(\sigma ,\tau)\right]_{\sigma ,\tau\preceq \alpha }$
is invertible. This implies that $\cM _K=0$ and $\cP ^0_N$
can be viewed as a subspace of $L^2(K)$. Also, for any 
$\alpha \in \FF _N^+$, $\{X_{\sigma }\}_{\sigma \preceq \alpha }$
is a linearly independent family in $L^2(K)$. Then, the 
Gram-Schmidt procedure gives a family $\{\varphi _{\sigma }\}
_{\sigma \in \FF _N^+}$ of elements in $\cP ^0_N$ such that 

\begin{equation}\label{cond1}
\varphi _{\sigma }=\sum _{\tau \preceq \sigma }a_{\sigma ,\tau}X_{\tau},
\quad a_{\sigma ,\sigma }>0;
\end{equation}
\begin{equation}\label{cond2}
\langle \varphi _{\sigma }, \varphi _{\tau }\rangle _K=0,
\quad \emptyset \preceq \sigma \prec \tau .
\end{equation}

An explicit formula for the orthogonal polynomials  
$\varphi _{\sigma }$ can be obtained in the same manner as in 
the classical (one variable) case.
Define for $\sigma \in \FF _N^+$, 
\begin{equation}\label{D}
D_{\sigma }=\det\left[K(\sigma ',\tau ')\right]_
{\sigma ',\tau '\preceq \sigma }
\end{equation} 
and let $\{\gamma _{\sigma }\}
_{\sigma \in \FF _N^+}$ be the parameters associated to $K$
as described in Section~2.2.3. Note that since 
all the matrices $\left[K(\sigma ,\tau)\right]_{\sigma ,\tau\preceq \alpha }$, 
$\alpha \in \FF _N^+$, are assumed to be invertible, it follows that  
$|\gamma _{\sigma }|<1$ for all $\sigma \in \FF _N^+$.

\begin{theorem}\label{T21}
$(1)$ $\varphi _{\emptyset }=1$ and for $\emptyset \prec \sigma $,
\begin{equation}\label{prop1}
\varphi _{\sigma }=\frac{1}{\sqrt{D_{\sigma -1}D_{\sigma }}}
{\det \left[\begin{array}{c}
\left[K(\sigma ',\tau ')\right]_
{\sigma '\prec \sigma ;
\tau '\preceq \sigma  } \\
  \\
\begin{array}{cccc}
1 & X_1 &  \ldots & X_{\sigma }
\end{array}
\end{array}
\right]}.
\end{equation}
$(2)$ For $\emptyset \prec \sigma =i_1\ldots i_k$,
\begin{equation}\label{pro2}
\varphi _{\sigma }=
\frac{1}{\prod _{1\leq j\leq k}(1-|\gamma _{i_j\ldots i_k}|^2)^{1/2}}
(X_{\sigma }
+\mbox{lower order terms}).
\end{equation}
\end{theorem}
\begin{proof}
The proof is similar to the classical one. Thus, we deduce 
from the orthogonality condition \eqref{cond2} that
$\langle \varphi _{\sigma },X_{\tau '}\rangle _K=0$
for $\emptyset \preceq \tau '\prec \sigma $, which implies 
that $\sum _{\tau \preceq \sigma }a_{\sigma ,\tau }
K(\tau ',\tau )=0$ for $\emptyset \preceq \tau '\prec \sigma $.
Using the Cramer rules for the system
$$\left\{\begin{array}{rcl}
\sum _{\tau \preceq \sigma }a_{\sigma ,\tau }
K(\tau ',\tau )&=&0, \quad \emptyset \preceq \tau '\prec \sigma ,\\
 & & \\
\sum _{\tau \preceq \sigma }a_{\sigma ,\tau }X_{\tau }&=&\varphi _{\sigma },
\end{array}
\right.
$$
with unknowns $a_{\sigma ,\tau }$, we deduce 
$$a_{\sigma ,\sigma }=
\frac{\varphi _{\sigma }D_{\sigma -1}}
{\det \left[\begin{array}{c}
\left[K(\sigma ',\tau ')\right]_
{\sigma '\prec \sigma ;\tau '\preceq \sigma } \\
 \\
\begin{array}{cccc}
1 & X_1 &  \ldots & X_{\sigma }
\end{array}
\end{array}
\right]}.$$
Therefore,
$$\varphi _{\sigma }=\frac{a_{\sigma ,\sigma }}{D_{\sigma -1}}
\det \left[\begin{array}{c}
\left[K(\sigma ',\tau ')\right]_
{\sigma '\prec \sigma ;\tau '\preceq \sigma } \\
\\
\begin{array}{cccc}
1 & X_1 & \ldots & X_{\sigma }
\end{array}
\end{array}
\right].$$
We now compute $a_{\sigma ,\sigma }$ and $D_{\sigma }$
in terms of the parameters $\{\gamma _{\sigma }\}_{\sigma \in \FF _N^+}$
of $K$.
First we notice that 
$$\langle \det \left[\begin{array}{c}
\left[K(\sigma ',\tau ')\right]_
{\sigma '\prec \sigma ;\tau '\preceq \sigma } \\
\\
\begin{array}{cccc}
1 & X_1 &  \ldots  & X_{\sigma }
\end{array}
\end{array}
\right],X_{\sigma }\rangle _K=D_{\sigma }$$
and since $X_{\sigma }=\frac{1}{a_{\sigma ,\sigma }}\varphi _{\sigma }+
\sum _{\tau \prec \sigma }c_{\tau }X_{\tau }$,
we deduce 
$$D_{\sigma }=\langle \frac{D_{\sigma -1}}{a_{\sigma ,\sigma }}
\varphi _{\sigma },\frac{1}{a_{\sigma ,\sigma }}\varphi _{\sigma }+
\sum _{\tau \prec \sigma }c_{\tau }X_{\tau }\rangle _K=\frac{D_{\sigma -1}}
{a^2_{\sigma ,\sigma }},$$
so that 
$$\frac{1}{a^2_{\sigma ,\sigma }}=\frac{D_{\sigma }}{D_{\sigma -1}}$$
which gives \eqref{prop1}. 

In order to compute $D_{\sigma }$ in terms of
$\{\gamma _{\sigma }\}_{\sigma \in \FF _N^+}$
we use Theorem~1.5.10 in \cite{Co} and the special
structure of $D_{\sigma }$. Thus, 
$$D_{\sigma }=\prod _{\emptyset \prec \sigma ',\tau '\preceq \sigma }
(1-|\gamma _{\sigma ',\tau '}|^2)$$
and for $\emptyset \prec \sigma =i_1\ldots i_k$ we deduce
$$\frac{1}{a^2_{\sigma ,\sigma }}=\frac{D_{\sigma }}{D_{\sigma -1}}=
\prod _{1\leq j\leq k}(1-|\gamma _{i_j\ldots i_k}|^2).$$
Then,
$$\begin{array}{rcl}
\varphi _{\sigma }&=&a_{\sigma ,\sigma }X_{\sigma }+ 
\sum _{\tau \prec \sigma }a_{\sigma ,\sigma }c_{\tau }X_{\tau }\\
 & & \\
 &=&\frac{1}{\prod _{1\leq j\leq k}(1-|\gamma _{i_j\ldots i_k}|^2)^{1/2}}
(X_{\sigma }
+\mbox{lower order terms}),
\end{array}$$
which gives \eqref{pro2}.
\end{proof}
We illustrate this 
result 
for $N=2$. From now on it is convenient to use the 
notation $d_{\sigma }=(1-|\gamma _{\sigma }|^2)^{1/2}$,
$\sigma \in \FF _N^+-\{\emptyset \}$.

\bigskip
\noindent
{\bf Example.} Let $N=2$ and assume the positive kernel $K$ satisfies the 
conditions in Theorem~\ref{T21}. We have $D_{\emptyset }=1$ and the next 
three determinants are:
$$D_1=\det \left[\begin{array}{cc}
1 & s_1 \\
\overline{s}_1 & 1
\end{array}\right]=d_1^2;
$$
$$D_2=\det \left[\begin{array}{ccc}
1 & s_1 & s_2 \\
\overline{s}_1 & 1 & 0 \\
\overline{s}_2 & 0 & 1
\end{array}\right]=d_1^2d_2^2;
$$
$$D_{11}=\det \left[\begin{array}{cccc}
1 & s_1 & s_2 & s_{11}\\
\overline{s}_1 & 1 & 0 & s_1 \\
\overline{s}_2 & 0 & 1 & 0 \\
\overline{s}_{11} & \overline{s}_1 & 0 & 1 
\end{array}\right]=d_1^4d_2^2d^2_{11}.
$$
Using Theorem~\ref{T21} we can easily calculate the first 
four orthogonal polynomials of $K$. Thus, $\varphi _{\emptyset }=1$
and then
$$\varphi _1=\frac{1}{d_1}\det 
\left[\begin{array}{cc}
1 & s_1 \\
1 & X_1
\end{array}\right]=-\frac{\gamma _1}{d_1}+\frac{1}{d_1}X_1;
$$
$$\varphi _2=\frac{1}{d^2_1d_2}\det 
\left[\begin{array}{ccc}
1 & s_1 & s_2 \\
\overline{s}_1 & 1 & 0 \\
1 & X_1 & X_2
\end{array}\right]=-\frac{\gamma _2}{d_1d_2}+
\frac{\overline{\gamma }_1\gamma _2}{d_1d_2}X_1+\frac{1}{d_2}X_2,
$$
where we used the fact that $s_2=d_1\gamma _2$. Then, after some 
calculations, 
$$\begin{array}{rcl}
\varphi _{11}&=&\displaystyle\frac{1}{d^3_1d^2_2d_{11}}\det 
\left[\begin{array}{cccc}
1 & s_1 & s_2 & s_{11} \\
\overline{s}_1 & 1 & 0 & s_1 \\
\overline{s}_2 & 0 & 1 & 0 \\
1 & X_1 & X_2 & X^2_1
\end{array}\right]\\
 & & \\
 &=&
-\displaystyle\frac{\gamma _{11}}{d_1d_2d_{11}}+
(-\displaystyle\frac{\gamma _1}{d_1d_{11}}+\displaystyle\frac{\gamma _{11}
\overline{\gamma }_1}{d_1d_2d_{11}})X_1+
\displaystyle\frac{\gamma _{11}\overline{\gamma }_2}{d_2d_{11}}X_2+
\displaystyle\frac{1}{d_{11}d_1}X_1^2.
\end{array}
$$

We establish that the orthogonal polynomials introduced above
satisfy equations similar to the classical Szeg\"o
difference equations. 

\begin{theorem}\label{T2}
The orthogonal polynomials satisfy the following recurrences:
$\varphi _{\emptyset }=1$ and for $k\in \{1,\ldots ,N\}$,
$\sigma \in \FF _N^+$, 
\begin{equation}\label{szego}
\varphi _{k\sigma }=\frac{1}{d_{k\sigma }}
(X_k\varphi _{\sigma }-\gamma _{k\sigma } 
\varphi ^{\sharp}_{k\sigma -1}),
\end{equation}
where $\varphi ^{\sharp}_{\emptyset }=1$ and for $k\in \{1,\ldots ,N\}$,
$\sigma \in \FF _N^+$,
\begin{equation}\label{sarp}
\varphi ^{\sharp}_{k\sigma }=\frac{1}{d_{k\sigma }}
(-\overline{\gamma }_{k\sigma }X_k\varphi _{\sigma }+ 
\varphi ^{\sharp}_{k\sigma -1}).
\end{equation}
\end{theorem}
\begin{proof}
We deduce this result from similar formulae obtained 
for an arbitrary positive definite kernel. 
In this way we can show the meaning of the polynomials 
$\varphi ^{\sharp}_{\sigma }$, $\sigma \in \FF _N^+$.
Let $\left[t_{i,j}\right]_{i,j\geq 1}$ be a positive definite kernel on 
$\NN $ and assume that each matrix
$A^{(i,j)}=\left[t_{k,l}\right]_{i\leq k,l\leq j}$, $1\leq i\leq j,$ 
is invertible. Also, assume $t_{k,k}=1$ for all $k\geq 1$.
Let $F_{i,j}$ be the upper Cholesky factor of $A^{(i,j)}$, so that
$F_{i,j}$ is an upper triangular matrix with positive diagonal and 
$A^{(i,j)}=F^*_{i,j}F_{i,j}$. A dual, lower Cholesky factor is obtained as
follows: define the symmetry of appropriate dimension,
$$\cJ =\left[\begin{array}{ccccc}
 0 & 0 & \ldots & 0 & I \\
 0 & 0 &        & I & 0 \\
 \vdots & & \ddots & & \\ 
 0 & I &        &  & 0 \\
 I & 0 &        & 0 & 0 
\end{array} 
\right]
$$ 
and then let $\tilde F_{i,j}$ denote the upper Cholesky factor of 
$B^{(i,j)}=\cJ A^{(i,j)}\cJ$. If $G_{i,j}=\cJ \tilde F_{i,j}\cJ$, then 
$$A^{(i,j)}=\cJ B^{(i,j)}\cJ=\cJ F^*_{i,j}F_{i,j}\cJ=G^*_{i,j}G_{i,j},$$
and $G_{i,j}$ is a lower triangular matrix with positive diagonal, 
called the lower Cholesky factor of $A^{(i,j)}$.
Let $P_{i,j}$ be the last column of $F^{-1}_{i,j}$ and let 
$P^{\sharp }_{i,j}$ be the first column of $G^{-1}_{i,j}$, that is
$$P_{i,j}=F^{-1}_{i,j}E, \quad 
P^{\sharp }_{i,j}=G^{-1}_{i,j}\cJ E,
$$
where $E=
\left[\begin{array}{cccc}
0 & \ldots & 0 & I
\end{array}
\right]^T$.
Let $\{r_{i,j}\}_{1<i\leq j}$ be the parameters associated to
$\left[t_{i,j}\right]_{i,j\geq 1}$ by Theorem~1.5.3 in \cite{Co}
and let $\rho _{i,j}=(1-|r_{i,j}|^2)^{1/2}$. 
We have that 
\begin{equation}\label{aa}
P_{1,n}=\frac{1}{d_{1,n}}\left[
\begin{array}{c}
0 \\
P_{2,n}
\end{array}
\right]-
\frac{r_{1,n}}{d_{1,n}}\left[
\begin{array}{c}
P^{\sharp }_{1,n-1} \\
0
\end{array}
\right],
\end{equation}
\begin{equation}\label{bb}
P^{\sharp }_{1,n}=-\frac{\overline{r}_{1,n}}{d_{1,n}}\left[
\begin{array}{c}
0 \\
P_{2,n}
\end{array}
\right]+
\frac{1}{d_{1,n}}\left[
\begin{array}{c}
P^{\sharp }_{1,n-1} \\
0
\end{array}
\right].
\end{equation}
These formulae are presumable known to the experts. 
For the sake of completeness we give a proof here based on results 
and notation from \cite{Co}.
First we introduce the following elements.
For $i<j$, 
\begin{equation}\label{linie}
L_i^{(j)}=L(\{r_{i,k}\}_{k=i+1}^j)=\left[
\begin{array}{ccccc}
r_{i,i+1} & \rho _{i,i+1}r_{i,i+1} & \ldots & 
\rho _{i,i+1}\ldots \rho _{i,j-1}r_{i,j}
\end{array}
\right];
\end{equation}
$$C_j^{(i)}=\left[
\begin{array}{c}
r_{j-1,j} \\
\vdots \\
r_{i+1,j}\rho _{i+2,j} \ldots \rho _{j-1,j} \\
r_{i,j}\rho _{i+1,j}\ldots \rho _{j-1,j}
\end{array}
\right],
$$
and 
$$K_i^{(j)}=\left[
\begin{array}{c}
\overline{r}_{i,i+1}\rho _{i,i+2}\ldots \rho _{i,j} \\
\vdots \\
\overline{r}_{i,j-1}\rho _{i,j} \\
\overline{r}_{i,j}
\end{array}
\right]=\left[\begin{array}{c}
K_i^{(j-1)}\rho _{i,j} \\
\overline{r}_{i,j}
\end{array}
\right].
$$
Also, we define inductively: $D_i^{(i+1)}=\rho _{i,i+1}$,
\begin{equation}\label{defect}
D_i^{(j)}=D(\{r_{i,k}\}_{k=i+1}^j)=\left[\begin{array}{cc}
D_i^{(j-1)} & -K_i^{(j-1)}r_{i,j} \\
0 & \rho _{i,j}
\end{array}
\right].
\end{equation}
We also need to review the factorization of unitary matrices.
This is an extension of Euler's description of $SO(3)$. First we 
define 
$$R_{j-i}(r_{i,k})=I_{k-1-i}\oplus \left[\begin{array}{cc}
r_{i,k} & \rho _{i,k} \\
\rho _{i,k} & -\overline{r}_{i,k}
\end{array}
\right]\oplus I_{j-k-1},
$$
where $I_{k-1-i}$ is the identity matrix of size $k-1-i$.
Then, 
$$R_{i,j}=R_{j-i}(r_{i,i+1})\ldots 
R_{j-i}(r_{i,j}),$$
and 
$$U_{i,j}=R_{i,j}(U_{i+1,j}\oplus 1).$$
It turns out that any unitary matrix can be written as a
matrix of the form of $U_{i,j}$. The main idea for the 
proof of \eqref{aa} is to use the identity
\begin{equation}\label{bafta}
U_{i,j}\cJ G_{i,j}=F_{i,j},
\end{equation}
which follows from the relations (1.6.10), (6.3.8) and 
(6.3.9) in \cite{Co}.
Thus, we notice that \eqref{bafta} implies  
$$P^{\sharp }_{i,j}=F_{i,j}U_{i,j}E,$$
which is more tractable than the original definition of $P^{\sharp }_{i,j}$.
This is seen from the following calculations. Using formula 
(1.5.7) in \cite{Co}, the above definition of $D_1^{(n)}$,
and the notation $D_1^{-(n)}=(D_1^{(n)})^{-1}$, we obtain
that 
$$\begin{array}{rcl}
P_{1,n}&=&
\left[\begin{array}{cc}
1 & -L_1^{(n)}D_1^{-(n)} \\
  & \\
0 & F_{2,n}^{-1}D_1^{-(n)}
\end{array}
\right]E=\left[\begin{array}{c}
 -L_1^{(n)}D_1^{-(n)}E \\
  \\
F_{2,n}^{-1}D_1^{-(n)}E
\end{array}
\right] \\
 & & \\ 
 &=&\left[\begin{array}{c}
 -L_1^{(n)}
\left[\begin{array}{cc}
D_1^{-(n-1)} & 
\frac{r_{1,n}}{\rho _{1,n}}
D_1^{-(n-1)}K_1^{(n-1)} \\
0 & \frac{1}{\rho _{1,n}}
\end{array}\right]E \\
  \\
F_{2,n}^{-1}
\left[\begin{array}{cc}
D_1^{-(n-1)} & 
\frac{r_{1,n}}{\rho _{1,n}}
D_1^{-(n-1)}K_1^{(n-1)} \\
0 & \frac{1}{\rho _{1,n}}
\end{array}\right]E
\end{array}
\right] \\
 & & \\ 
 &=&\left[\begin{array}{c}
 -L_1^{(n)}
\left[\begin{array}{c}
\frac{r_{1,n}}{\rho _{1,n}}
D_1^{-(n-1)}K_1^{(n-1)} \\
\frac{1}{\rho _{1,n}}
\end{array}\right] \\
 \\
F_{2,n}^{-1}
\left[\begin{array}{c}
\frac{r_{1,n}}{\rho _{1,n}}
D_1^{-(n-1)}K_1^{(n-1)} \\
\frac{1}{\rho _{1,n}}
\end{array}\right]
\end{array}
\right] \\
 & & \\
 &=&\displaystyle\frac{1}{\rho _{1,n}}
\left[\begin{array}{c}
0 \\
 \\
F_{2,n}^{-1}E
\end{array}
\right]
+
\left[\begin{array}{c}
 -L_1^{(n)}
\left[\begin{array}{c}
\frac{r_{1,n}}{\rho _{1,n}}D_1^{-(n-1)}K_1^{(n-1)} \\
\frac{1}{\rho _{1,n}}
\end{array}\right] \\
 \\
F_{2,n}^{-1}
\left[\begin{array}{c}
\frac{r_{1,n}}{\rho _{1,n}}{D_1^{-(n-1)}}^{-1}K_1^{(n-1)} \\
 0
\end{array}\right]
\end{array}
\right] \\
 &=&
\displaystyle\frac{1}{\rho _{1,n}}
\left[\begin{array}{c}
0 \\
 \\
P_{2,n}
\end{array}
\right]
+
\displaystyle\frac{r_{1,n}}{\rho _{1,n}}
\left[\begin{array}{c}
 -L_1^{(n-1)}D_1^{-(n-1)}K_1^{(n-1)}
-\rho _{1,2}\ldots \rho _{1,n-1} \\
\\
F_{2,n}^{-1}
\left[\begin{array}{c}
D_1^{-(n-1)}K_1^{(n-1)} \\
 0
\end{array}\right]
\end{array}
\right]. 
\end{array}
$$
The proof of formula (1.6.15) in \cite{Co} gives
$$L_1^{(n-1)}D_1^{-(n-1)}K_1^{(n-1)}
+\rho _{1,2}\ldots \rho _{1,n-1}
=\frac{1}{\rho _{1,2}\ldots \rho _{1,n-1}}
$$
and using formula (1.5.6) in \cite{Co} we deduce
$$F_{2,n}^{-1}
\left[\begin{array}{c}
D_1^{-(n-1)}K_1^{(n-1)} 
\\
 0
\end{array}\right]=
\left[\begin{array}{c}
F_{2,n-1}^{-1}D_1^{-(n-1)}K_1^{(n-1)} \\
0 
\end{array}\right],
$$
therefore
$$P_{1,n}=
\displaystyle\frac{1}{\rho _{1,n}}
\left[\begin{array}{c}
0 \\
P_{2,n}
\end{array}
\right]
-
\displaystyle\frac{r_{1,n}}{\rho _{1,n}}
\left[\begin{array}{c}
\frac{1}{\rho _{1,2}\ldots \rho _{1,n-1}} \\
\\
-F_{2,n-1}^{-1}D_1^{-(n-1)}K_1^{(n-1)} \\
 \\
0 
\end{array}\right].$$

\noindent
It remains to show that 
$$P^{\sharp }_{1,n-1}=
\left[\begin{array}{c}
\frac{1}{\rho _{1,2}\ldots \rho _{1,n-1}} \\
\\
-F_{2,n-1}^{-1}D_1^{-(n-1)}K_1^{(n-1)} 
\end{array}\right].$$
To that end we notice that using formula (1.5.8) in \cite{Co},
the definition of $U_{1,n-1}$, the fact that $R_{1,n-1}$
is a unitary matrix, and the notation 
$L_1^{*(n-1)}=(L_1^{(n-1)})^*$, we obtain 
$$\begin{array}{rcl}
P^{\sharp }_{1,n-1}&=&
\left[\begin{array}{cc}
0 & \frac{1}{\rho _{1,2}\ldots \rho _{1,n-1}} \\
 & \\
F_{2,n-1}^{-1} & -\frac{1}{\rho _{1,2}\ldots \rho _{1,n-1}}
F_{2,n-1}^{-1}L_1^{*(n-1)}
\end{array}\right]R^*_{1,n-1}R_{1,n-1}
\left[
\begin{array}{cc}
U_{2,n-1} & 0 \\
0 & 1 
\end{array}\right]E \\
 & & \\
 &=& 
\left[\begin{array}{c}
\frac{1}{\rho _{1,2}\ldots \rho _{1,n-1}} \\
\\
-\frac{1}{\rho _{1,2}\ldots \rho _{1,n-1}}
F_{2,n-1}^{-1}L_1^{*(n-1)})
\end{array}\right].
\end{array}
$$

\noindent
It follows that all we have to show is the equality
$$
F_{2,n-1}^{-1}D_1^{-(n-1)}K_1^{(n-1)}=
\frac{1}{\rho _{1,2}\ldots \rho _{1,n-1}}
F_{2,n-1}^{-1}L_1^{*(n-1)}.
$$
Now this is a simple consequence of the formula
$T^*D_{T^*}=D_{T}T^*$ for the contraction
$T=L_1^{(n-1)}$.
The formula \eqref{bb} can be proved in a similar manner.

We rewrite \eqref{aa} and \eqref{bb} for a positive definite
kernel $K$ satisfying \eqref{sta1} and \eqref{sta2}. We notice
that $P^{\sharp }_{1,n-1}$ is replaced by 
$\varphi ^{\sharp }_{k\sigma -1}$ and then we have to show that 
$P_{2,n}$ can be expressed in terms of $\varphi _{\sigma }$.
This follows by taking into account the relations  
\eqref{sta3} and \eqref{sta4} and using systematically 
\eqref{aa}. We can omit the details.
\end{proof}

The previous recurrence equations look quite similar to
the classical Szeg\"o recursions, still they have a 
new component that is best understood when 
considering \eqref{aa} and \eqref{bb}. This type of
recurrence equations was also found in \cite{CJ1}
in connection with some derivations on $\cL _N$. It turns out
that these derivations are related to those considered
in \cite{Ha}, and later studied in \cite{Fo}.

We also notice a graded form of the recurences \eqref{szego} and \eqref{sarp}.
It is convenient to introduce the following notation: for $n\geq 1$,
we use \eqref{linie} and define 
$$g_n=L(\{\gamma _{\sigma }\}_{|\sigma |=n}).$$
It was explained in \cite{CJ2} that $g_n$, $n\geq 1$, are 
the parameters associated to the kernel $K$ in \cite{Po}.
We also use \eqref{defect} in order to introduce the notation
$$H_n=D(\{\gamma _{\sigma }\}_{|\sigma |=n}), \quad n\geq 1.$$ 
Let $\sigma (n)$
be the largest word (with respect to the lexicographic order)
of lenght $n$, that is 
$\sigma (n)=\underbrace{N\ldots N}_
{N \,terms}$.

\begin{corollary}\label{graded}
The Szeg\"o polynomials satisfy the recurrences: for $k\geq 1$,
\begin{equation}\label{gradedszego}
\left[\varphi _{|\sigma |=k}\right]=
\left(
\left[\begin{array}{ccc}
X_1 & \ldots & X_N
\end{array}
\right] 
\left[\varphi _{|\sigma |=k-1}\right]^{\oplus N}-
\varphi ^{\sharp }_{\sigma (k)-1}g_k\right)H_k^{-1};
\end{equation}
\begin{equation}\label{gradedsarp}
\varphi ^{\sharp }_{\sigma (k)}=
\prod _{|\tau |=k}d_{\tau }^{-1}\left(
-\left[\begin{array}{ccc}
X_1 & \ldots & X_N
\end{array}
\right] 
\left[\varphi _{|\sigma |=k-1}\right]^{\oplus N}g^*_k+
\varphi ^{\sharp }_{\sigma (k)-1}\right).
\end{equation}
\end{corollary}
\begin{proof}
Both \eqref{gradedszego} and \eqref{gradedsarp}
follow by direct calculations from Theorem~\ref{T2}. 
\end{proof}

\section{Christoffel-Darboux formula}
A first consequence of the Szeg\"o formula in the classical case
is the Christoffel-Darboux formula. Here we find a similar
formula in several noncommuting variables. 
To that end we introduce additional notation. Let $\cE $ be 
a Hilbert space. In this paper $\cE $ will always be infinite 
dimensional. The 
$N$-dimensional unit ball of $\cE $
is defined by
$$\cB _N(\cE )=\{Z=\left(\begin{array}{ccc}
Z_1 & \ldots & Z_N
\end{array}
\right)\mid (Z\mid Z)<I_{\cE}\},
$$
where for two elements $Z=\left(\begin{array}{ccc}
Z_1 & \ldots & Z_N
\end{array}
\right)$ and 
$W=\left(\begin{array}{ccc}
W_1 & \ldots & W_N
\end{array}
\right)$
in $\cL (\cE )^N$ we define
\begin{equation}\label{inner}
(Z\mid W)=\sum _{k=1}^NZ_kW^*_k.
\end{equation}
We also need a sort of Szeg\"o kernel for $\cB _N(\cE )$.
One suggestion was given in \cite{CJ1} to consider
the following construction. For $Z\in \cB _N(\cE )$
define 
\begin{equation}\label{E}
E(Z)=\left[Z_{\sigma }\right]_{|\sigma |=0}^{\infty }
\in \cL (\oplus _{k\geq 0}\cE _k,\cE ).
\end{equation}
Also we use the notation $\mbox{diag}(S)$ to denote
the diagonal operator in $\cL (\oplus _{k\geq 0}\cE _k)$
with diagonal $S$. A Szeg\"o type kernel on $\cB _N(\cE )$
is given by the formula 
$$K_S(Z,W)=E(Z)E(W)^*,\quad Z,W\in \cB _N(\cE ).$$
The next result explain two important properties of $K$.

\begin{lemma}\label{duminica}
$(a)$ For any $T\in \cL (\cE )$ and $Z,W\in \cB _N(\cE )$, 
$$E(Z)\mbox{diag}(T-\sum _{k=1}^NZ_kTW^*_k)E(W)^*=T.$$
\smallskip
\noindent
$(b)$ The set $\{E(W)^*\cE \mid W\in \cB _N(\cE )\}$
is total in $\oplus _{k\geq 0}\cE _k$.
\end{lemma}
\begin{proof} $(a)$  Using directly the definitions,
$$\begin{array}{l}
 E(Z)\mbox{diag}(T-\sum _{k=1}^NZ_kTW^*_k)E(W)^* \\
  \\
 =T+\sum _{|\sigma |\geq 1}Z_{\sigma }TW^*_{\sigma }-
\sum _{k=1}^NE(Z))\mbox{diag}(Z_kTW^*_k)E(W)^* \\
  \\
 =T+\sum _{|\sigma |\geq 1}Z_{\sigma }TW^*_{\sigma }-
\sum _{k=1}^N\sum _{|\sigma |\geq 0}Z_{\sigma }Z_kTW^*_k
W^*_{\sigma } \\
  \\
 =T+\sum _{|\sigma |\geq 1}Z_{\sigma }TW^*_{\sigma }-
\sum _{|\sigma |\geq 1}Z_{\sigma }TW^*_{\sigma } \\
  \\
 =T.
\end{array}
$$

$(b)$ Let $e=\{e_{\sigma }\}_{\sigma \in \FF ^+_N}$
be an element of $\oplus _{k\geq 0}\cE _k$ orthogonal
to the linear span of $\{E(W)^*\cE \mid W\in \cB _N(\cE )\}$.
Taking $W=0$, we deduce that $e_{\emptyset }=0$. Next, we claim
that for each $\sigma \in \FF ^+_N-\{\emptyset \}$ there exist
$$W_l=(W_1^l,\ldots ,W_N^l)\in \cB _N(\cE ), \quad 
l=1,\ldots ,2|\sigma |,$$
such that 
$$\mbox{range}\left[
\begin{array}{ccc}
W_{\sigma }^{*1} & \ldots & W_{\sigma }^{*2|\sigma |}
\end{array}\right]=\cE,$$
and 
$$W_{\tau }^l=0\quad \mbox{for all}\quad 
\tau \ne \sigma , \quad |\tau |\geq |\sigma |, \quad l=1,\ldots ,2|\sigma |.$$

Once this claim is proved, a simple inductive argument gives
$e=0$, so 
$\{E(W)^*\cE \mid W\in \cB _N(\cE )\}$
is total in $\oplus _{k\geq 0}\cE _k$.
Therefore we focus on the proof of the claim.

Let $\{e^n_{ij}\}_{i,j=1}^n$ be the matrix units of the algebra
$M_n$ of $n\times n$ matrices. Each 
$e^n_{ij}$ is an $n\times n$ matrix consisting of $1$
in the $(i,j)th$ entry and zeros elsewhere.
For a Hilbert space $\cE _1$ we define $E^n_{ij}=e^n_{ij}\otimes 
I_{\cE _1}$ and we notice that 
\begin{equation}\label{units}
E^n_{ij}E^n_{kl}=\delta _{jk}E^n_{il},\quad E^{*n}_{ji}=E^n_{ij}.
\end{equation}

Let $\cE $ be infinite dimensional, so that $\cE =\cE _1^{\oplus 2|\sigma |}$
for some Hilbert space $\cE _1$. Let $\sigma =i_1\ldots i_k$.
For $s=1,\ldots ,N$, we define 
$$J_s=\{l\in \{1,\ldots ,k\}\mid i_{k+1-l}=s\}$$
and 
$$
W^{*p}_{s }=\frac{1}{\sqrt{2}}
\sum _{r\in J_s}E^{2|\sigma |}_{r+p-1,r+p},\quad 
s=1,\ldots ,N, \quad p=1,\ldots ,|\sigma |.$$
We show that for each $p\in \{1,\ldots ,|\sigma |\}$,
\begin{equation}\label{kiwi}
W^{*p}_{\sigma }=\frac{1}{{\sqrt{2^k}}}E^{2|\sigma |}_{p,k+p},
\end{equation}
\begin{equation}\label{mango}
W_{\tau }^p=0\quad \mbox{for} \quad 
\tau \ne \sigma ,\quad |\tau |\geq |\sigma |.
\end{equation}
Using
\eqref{units}, we deduce 
$$\begin{array}{rcl}
\sum _{s=1}^NW_s^pW_s^{*p}&=&
\frac{1}{2}\sum _{s=1}^N\sum _{r\in J_s}E^{2|\sigma |}_{r+p,r+p-1}
E^{2|\sigma |}_{r+p-1,r+p} \\
 & & \\
 &=&\frac{1}{2}\sum _{s=1}^N\sum _{r\in J_s}E^{2|\sigma |}_{r+p,r+p} \\
 & & \\
 &=&\frac{1}{2}\sum _{r=1}^kE^{2|\sigma |}_{r+p,r+p}<I,
\end{array}
$$
hence $W^p\in \cB _N(\cE )$ for each $p=1,\ldots ,|\sigma |$.
For each word $\tau =j_1\ldots j_k\in \FF _N^+-\{\emptyset \}$
we deduce by induction that 
\begin{equation}\label{apple}
W^{*p}_{j_k}\ldots 
W^{*p}_{j_1}=\frac{1}{\sqrt{2^k}}
\sum _{r\in A_{\tau }}E^{2|\sigma |}_{r+p-1,r+p+k-1},
\end{equation}
where 
$A_{\tau }=\cap _{p=0}^{k-1}(J_{j_{k-p}}-p)\subset \{1,\ldots ,N\}$
and $J_{j_{k-p}}-p=\{l-p\mid l\in J_{i_{k-p}}\}$.

We show that $A_{\sigma }=\{1\}$ and $A_{\tau }=\emptyset $
for $\tau \ne \sigma $. Let $q\in A_{\tau }$. Therefore, for any 
$p\in \{0,\ldots ,k-1\}$ we must have $q+p\in J_{j_{k-p}}$
or $i_{k+1-q-p}=j_{k-p}$. For $p=k-1$ we deduce 
$j_1=i_{2-q}$ and since $2-q\geq 1$, it follows that $q\leq 1$. 
Also $q\geq 1$, therefore the only element that 
can be in $A_{\tau }$ is $q=1$, in which case we must have
$\tau =\sigma $. Since $l\in J_{i_{k+1-l}}$
for each $l=1, \ldots ,k-1$, we deduce that $A_{\sigma }=\{1\}$
and $A_{\tau }=\emptyset $ for $\tau \ne \sigma $.
Formula \eqref{apple} implies \eqref{kiwi}. In a similar manner
we can construct a family $W^p$, $p=|\sigma |+1, \ldots ,2|\sigma |,$
such that 
$$W^{*p}_{\sigma }=\frac{1}{\sqrt{2^k}}E^{2|\sigma |}_{p+k,p},$$
and 
$$
W_{\tau }^p=0\quad \mbox{for} \quad 
\tau \ne \sigma ,\quad |\tau |\geq |\sigma |.
$$
Thus, for $s=1,\ldots ,N$, we define 
$$K_s=\{l\in \{1,\ldots ,k\}\mid i_k=s\}$$
and 
$$W^{*p}_{s}=\frac{1}{\sqrt{2}}
\sum _{r\in K_s}E^{2|\sigma |}_{r+p-k,r+p-k-1},\quad 
s=1,\ldots ,N, \quad p=|\sigma |+1,\ldots ,2|\sigma |.$$

Now, 
$$\left[
\begin{array}{ccc}
W_{\sigma }^{*1} & \ldots & 
W_{\sigma }^{*2|\sigma |}
\end{array}
\right]=
\frac{1}{\sqrt{2^k}}\left[
\begin{array}{cccccc}
E^{2|\sigma |}_{1,k+1} & \ldots & 
E^{2|\sigma |}_{k,2k} & E^{2|\sigma |}_{k+1,1}
& \ldots & E^{2|\sigma |}_{2k,k}
\end{array}
\right],
$$
whose range is $\cE $. 
This concludes the proof. 
\end{proof}

We note that the result given by Lemma~\ref{duminica}(b)
is not true in case $\cE $ is finite dimensional. The meaning 
of the result is that in case $\cE $ is infinite dimensional
then $E$ is really the Kolmogorov decomposition of 
the kernel $K_S$.

We now let a polynomial $P=\sum _{\sigma \in \FF _N^+}c_{\sigma }
X_{\sigma }\in \cP ^0_N$ take values on $\cB _N(\cE )$ 
by the formula 
\begin{equation}\label{eval}
P(Z)=\sum _{\sigma \in \FF _N^+}c_{\sigma }
Z_{\sigma }, \quad Z\in \cB _N(\cE ).
\end{equation}

Define the Cristoffel-Darboux kernel by the formula
\begin{equation}\label{CD}
K_{CD}(Z,W)=
E(Z)\mbox{diag}(\varphi ^{\sharp }_{\sigma (n)}(Z)
\varphi ^{\sharp }_{\sigma (n)}(W)^*-\sum _{|\tau |=n}
\varphi _{\tau }(Z)
\varphi _{\tau }(W)^*)E(W)^*,
\end{equation}
for 
$Z,W\in \cB _N(\cE )$.

\begin{theorem}\label{cristica}
For any $Z,W\in \cB _N(\cE )$,
$$K_{CD}(Z,W)=\sum _{0\leq |\tau |<n}
\varphi _{\tau }(Z)
\varphi _{\tau }(W)^*.
$$
\end{theorem}
\begin{proof}
From \eqref{szego} and \eqref{sarp}
we deduce 
$$
\varphi ^{\sharp }_{k\sigma }(Z)
\varphi ^{\sharp }_{k\sigma }(W)^*-
\varphi _{k\sigma }(Z)
\varphi _{k\sigma }(W)^* 
=\varphi ^{\sharp }_{k\sigma -1}(Z)
\varphi ^{\sharp }_{k\sigma -1}(W)^*-
Z_k\varphi _{k\sigma }(Z)
\varphi _{k\sigma }(W)^*W^*_k,
$$
for any $k\in \{1,\ldots ,N\}$, $\sigma \in \FF ^+_N$, 
and $Z,W\in \cB _N(\cE )$.
Adding all these relations for  
$k\in \{1,\ldots ,N\}$ and $0\leq |\sigma |\leq n-1$, we deduce
$$
\varphi ^{\sharp }_{\sigma (n)}(Z)
\varphi ^{\sharp }_{\sigma (n)}(W)^*
-\sum _{0\leq |\sigma |\leq n}\varphi _{\sigma }(Z)
\varphi _{\sigma }(W)^* 
=\sum _{k=1}^N\sum _{0\leq |\sigma |\leq n-1}
Z_k\varphi _{\sigma }(Z)
\varphi _{\sigma }(W)^*W^*_k.
$$
This relation and Lemma \ref{duminica}
give 
$$\begin{array}{l}
     K_{CD}(Z,W) \\ 
     \\
     =E(Z)(\sum _{0\leq |\sigma |<n}\varphi _{\sigma }(Z)
\varphi _{\sigma }(W)^*-
\sum _{k=1}^NZ_k(\sum _{0\leq |\sigma |<n}\varphi _{\sigma }(Z)
\varphi _{\sigma }(W)^*)W^*_k)E(W)^* \\
  \\
=\sum _{0\leq |\tau |<n}
\varphi _{\tau }(Z)
\varphi _{\tau }(W)^*.
\end{array}
$$
\end{proof}

We can show one more application of
Lemma~\ref{duminica}. For a formal power series
$$f=\sum _{\sigma \in \FF ^+_N}c_{\sigma }X_{\sigma },$$
in $N$ noncommuting variables $X_1$, $\ldots $, $X_N$, 
we denote by $T_f$ the lower triangular infinite matrix
associated to $f$ as described in Section~2.1.  
We denote by $\cS _N$ the Schur class 
of those formal power series $f$ with the property that 
$T_f$ is a contraction in $\cL (\oplus _{k\geq 0}\CC _k)$.
If $\cE $ is an infinite dimensional
Hilbert space then we can define 
$f(Z)$ for $Z\in \cB _N(\cE)$ as in \cite{CJ1}, by the formula
\begin{equation}\label{evaluare}
f(Z)=E(Z)(T_f\otimes I_{\cE })/\cE .
\end{equation}
We notice that this definition is consistent with \eqref{eval}.
We extend a familiar characterization 
of the Schur class to the setting of this paper.
\begin{theorem}\label{cara}
The formal power series $f$ belongs to $\cS _N$
if and only if 
$$C_f(Z,W)=E(Z)\mbox{diag}(I-f(Z)f(W)^*)E(W)^*,\quad Z,W\in 
\cB _N(\cE),$$
is a positive definite kernel on $\cB _N(\cE)$.
\end{theorem}
\begin{proof}
Using Lemma~3.1 in \cite{CJ1} we deduce
that for $Z,W\in 
\cB _N(\cE)$,
$$\begin{array}{l}
  E(Z)(I-(T_f\otimes I_{\cE })(T_f\otimes I_{\cE })^*)E(W)^* \\
  \\
 =E(Z)E(W)^*-E(Z)(T_f\otimes I_{\cE })(T_f\otimes I_{\cE })^*)E(W)^* \\
  \\ 
 =E(Z)E(W)^*-E(Z)\mbox{diag}(f(Z))\mbox{diag}(f(W)^*)E(W)^* \\
  \\
 =E(Z)\mbox{diag}(I-f(Z)f(W)^*)E(W)^*=C_f(Z,W).
\end{array}
$$
This relation implies that if $f\in \cS _N$ then $C_f$
is a positive definite kernel on $\cB _N(\cE)$.
For the converse implication 
we have to use in addition Lemma~\ref{duminica}.
\end{proof}

\section{Inverse problems}
In this brief section we prove a Favard type 
result for orthogonal polynomials in several 
noncommuting variables.
\begin{theorem}
Let $\{\gamma _{\sigma }\}_{\sigma \in \FF ^+_N}$
be a family of complex numbers with $\gamma _{\emptyset }=0$
and $|\gamma _{\sigma }|<1$ for $\sigma \in \FF ^+_N-\{\emptyset \}$.
Then there exists a unique positive definite kernel
$K$ satisfying \eqref{sta1} and \eqref{sta2} such that the 
polynomials $\varphi _{\sigma }$, $\sigma \in 
\FF ^+_N$,  defined by the recursions: 
$\varphi _{\emptyset }=1$, $\varphi ^{\sharp}_{\emptyset }=1$,
and for $k\in \{1,\ldots ,N\}$,
$\sigma \in \FF _N^+$, 
$$
\varphi _{k\sigma }=\frac{1}{d_{k\sigma }}
(X_k\varphi _{\sigma }-\gamma _{k\sigma } 
\varphi ^{\sharp}_{k\sigma -1}),
$$
$$
\varphi ^{\sharp}_{k\sigma }=\frac{1}{d_{k\sigma }}
(-\overline{\gamma }_{k\sigma }X_k\varphi _{\sigma }+ 
\varphi ^{\sharp}_{k\sigma -1}),
$$
are orthogonal with respect to $K$.
\end{theorem}
\begin{proof}
Once again we rely on some results that are known for arbitrary
positive definite kernels on the set of integers. In this way, the 
proof is quite straightforward. Let $\{\gamma _{\sigma ,\tau }
\mid \sigma ,\tau \in \FF ^+_N, \sigma \preceq \tau \}$ be the 
family of complex numbers associated to 
$\{\gamma _{\sigma }\}_{\sigma \in \FF ^+_N}$ by \eqref{explo}.
Let $K$ be the positive definite kernel associated to 
$\{\gamma _{\sigma ,\tau }
\mid \sigma ,\tau \in \FF ^+_N, \sigma \preceq \tau \}$
by Theorem~1.5.3 in \cite{Co}. By Theorem~\ref{T2}, the polynomials 
$\varphi _{\sigma }$, $\sigma \in \FF ^+_N$, defined by the recurrences:
$\varphi _{\emptyset }=1$, $\varphi ^{\sharp}_{\emptyset }=1$,
and for $k\in \{1,\ldots ,N\}$,
$\sigma \in \FF _N^+$, 
$$
\varphi _{k\sigma }=\frac{1}{d_{k\sigma }}
(X_k\varphi _{\sigma }-\gamma _{k\sigma } 
\varphi ^{\sharp}_{k\sigma -1}),
$$
$$
\varphi ^{\sharp}_{k\sigma }=\frac{1}{d_{k\sigma }}
(-\overline{\gamma }_{k\sigma }X_k\varphi _{\sigma }+ 
\varphi ^{\sharp}_{k\sigma -1}),
$$
must be the orthogonal polynomials of $K$.
\end{proof}

\end{document}